\documentclass[a4paper,11pt,reqno]{article}
\usepackage{amsthm}
\usepackage{amsmath}
\usepackage{amssymb}
\usepackage[mathscr]{eucal}
\usepackage[latin2]{inputenc}
\usepackage[T1]{fontenc}
\usepackage{indentfirst}
\usepackage{amsfonts}
\usepackage{graphicx}
\usepackage{color}
\usepackage[top=3cm, bottom=3cm, left=3cm, right=3cm]{geometry}

\newtheoremstyle{new_plain}
	{}
	{}
	{\itshape}
	{}
	{\sffamily\bfseries}
	{.}
	{5pt plus1pt minus1pt\relax}
	{\thmnumber{#2. }\thmname{#1}\thmnote{ #3}}

\newtheoremstyle{nev_circle_definition}
	{}
	{}
	{\normalfont}
	{}
	{\sffamily\bfseries}
	{.}
	{5pt plus1pt minus1pt\relax}
	{\thmnumber{#2. }\thmname{#1}\thmnote{ #3}}
	
\theoremstyle{plain}
	\newtheorem{theorem}{Theorem}
	\newtheorem{lemma}[theorem]{Lemma}
	
	\newtheorem{remark}[theorem]{Remark}

\theoremstyle{definition}
	\newtheorem{definition}[theorem]{Definition}

\theoremstyle{nev_circle_definition}

\numberwithin{theorem}{section}
\numberwithin{equation}{section}

\title{\bfseries \Large Christoffel functions on Jordan curves with respect to measures with jump singularity}
\author{\Large Tivadar Danka\thanks{This research was supported by ERC Advanced Grant No. 267055} \\[0.3cm] \large Bolyai Institute, University of Szeged \\ Aradi V. tere 1, H-6720 Szeged, Hungary \\ email: tdanka@math.u-szeged.hu}
\date{}

\begin{document}

\maketitle

\begin{abstract}
	In this paper we establish asymptotic results for Christoffel functions with respect to measures supported on Jordan curves having a Radon-Nikodym derivative with a jump singularity. We extend the results known for measures supported in the interval \( [-1,1] \) to more general sets, for example a system of Jordan curves, using polynomial inverse images. It is shown that the asymptotic limit can be written in terms of the equilibrium measure or Green's function.
\end{abstract}


\section{Introduction}

Let \( \mu \) be a Borel measure with compact support in the complex plane, and assume that \( \operatorname{supp}(\mu) \) is an infinite set. The functions
\begin{equation}\label{christoffel_function}
	\lambda_n(\mu,z) = \inf_{P_n(z) = 1} \int |P_n|^2 d\mu,
\end{equation}
where the infimum is taken over all polynomials of degree at most \( n \), are called the Christoffel functions associated with \( \mu \). Let \( p_k(z) \) denote the orthonormal polynomials with respect to \( \mu \). It is known that
\begin{equation}\label{kernel_diagonal_sum}
	\lambda_{n}^{-1}(\mu,z) = \sum_{k=0}^{n} |p_k(z)|^2.
\end{equation}
Our goal is to study the asymptotic behavior of the Christoffel functions for a measure $ \mu $ supported on a finite union of Jordan curves at a point $ z_0 $ for which $ \mu $ is absolutely continuous in a subarc $ J $ containing $ z_0$ with respect to the arc length measure, and its density function $ w_\mu $  having a jump singularity, i.e. if \( f_J:(-1,1) \to J \) is the parametrization of $ J $ with \( f_J(0) = z_0 \), then
\begin{equation}\label{type_one_singularity}
	\lim_{t \to 0-} w_\mu(f_J(t)) = A, \quad \lim_{t \to 0+} w_\mu(f_J(t)) = B,
\end{equation}
for some \( A, B > 0 \), \( A \neq B \).

Christoffel functions have been studied extensively even in the last century. One important early due to Szeg\H{o} states that if \( \mu \) is an absolutely continuous measure on the unit circle with \( d\mu(t) = \mu'(t)dt \), then if \( \log \mu'(t) \) is integrable, we have
\[
	\lim_{n \to \infty} \lambda_n(\mu,z) = (1-|z|^2) \exp \Bigg( \frac{1}{2\pi} \int_{-\pi}^{\pi} \frac{e^{it} - z}{e^{it} + z} \log \mu'(t) dt \Bigg), \quad |z| < 1.
\]
For more details on this result, see the book \cite{GSz} or Section 4.11 in \cite{PN}.

Asymptotic results for Christoffel functions with respect to measures supported in the complex plane with continuous density function have been established in \cite{T1} by V. Totik. In \cite{MMS}, A. Foulqui\'e Moreno, A. Martinez-Finkelshtein and V. L. Souza proved asymptotic results for measures supported on \( [-1,1] \) where the density function has a jump singularity at \( 0 \). That result has been extended to measures supported on a compact subset of the real line by P. Nevai and V. Totik in \( \cite{NeT} \). Our aim is to establish similar results for measures with support on a finite union of rectifiable Jordan curves. In this section we recall a few important concepts and state our main theorem. We begin with the definitions of capacity and the equilibrium measure.

Let \( \mu \) be a finite Borel measure with compact support, and define the so-called energy functional
\[
	I(\mu) = \int \int \log |z-w| d\mu(z) d\mu(w).
\]
For a compact set \( K \subseteq \mathbb{C} \), the energy of \( K \) is defined as
\[
	I(K) = \sup\{I(\mu): \operatorname{supp}(\mu) = K, \mu(K) = 1 \}.
\]
The capacity of \( K \) is defined as
\[
	\operatorname{cap}(K) = e^{I(K)}.
\]
If \( E \subseteq \mathbb{C} \) is an arbitrary Borel set, then the capacity is defined as
\[
	\operatorname{cap}(E) = \sup\{ \operatorname{cap}(K): K \subseteq \mathbb{E}, K \textnormal{ compact} \}.
\]
For a detailed discussion on capacity, see \cite{TR} or \cite{ST}.
\begin{definition}
	Let \( K \) be a compact subset of \( \mathbb{C} \). If \( \omega_K \) is a measure such that \( I(\omega_K) = I(K) \), then \( \omega_K \) is called an equilibrium measure.
\end{definition}
For compact subsets \( K \subseteq \mathbb{C} \) of finite energy \( I(K) > -\infty \), the equilibrium measure exists and is unique. For a detailed discussion on equilibrium measures and its properties, see \( \cite{TR} \).
\begin{definition}
	Let \( \mu \) be a Borel measure such that its support \( K=\operatorname{supp}(\mu) \) is compact and contains infinitely many points. \( \mu \) belongs to the class \( \mathbf{Reg} \), if for any sequence \( \{ P_n \}_{n \in \mathbb{N}} \) of nonzero polynomials with \( \deg(P_n) \leq n \), we have
\[
	\limsup_{n \to \infty} \Bigg( \frac{|P_n(z)|}{\| P_n \|_{L^2(\mu)}} \Bigg)^{1/n} \leq 1, \quad z \in K \setminus E,
\]
where \( E \) is a set with capacity \( 0 \).
\end{definition}
There are multiple equivalent definitions for the \( \mathbf{Reg} \) class. For more details, see \( \cite{STOP} \), Theorems 3.1.1, 3.2.1 and 3.2.3. Most important for us is the case when the support is regular with respect to the Dirichlet problem.
\begin{remark}\label{Reg_Dirichlet}
	\textnormal{Let \( \mu \in \mathbf{Reg} \). Assume that \( K = \operatorname{supp}(\mu) \) and \( \Omega \) denotes the unbounded component of \( \mathbb{C} \setminus K \). If \( \mathbb{C} \setminus \Omega \) is regular with respect to the Dirichlet problem, then for every sequence \( \{ P_n \}_{n \in \mathbb{N}} \) of nonzero polynomials with \( \deg(P_n) \leq n \), we have the uniform estimate}
\begin{equation}\label{Reg_Dirichlet_estimate}
	\limsup_{n \to \infty} \Bigg( \frac{\| P_n \|_{K}}{\| P_n \|_{L^2(\mu)}} \Bigg)^{1/n} \leq 1.
\end{equation}
\end{remark}
In the proof of the main theorem, it is enough to use the uniform estimate (\ref{Reg_Dirichlet_estimate}). An important regularity criteria can be given with the aid of the Radon-Nikodym derivative with respect to the equilibrium measure. If \( \mu \) is a Borel measure on the complex plane with compact support \( \operatorname{supp}(\mu) = K \), then
\begin{equation}\label{Erdos-Turan}
	\frac{d\mu}{d\omega_K} > 0 \quad \omega_K -a.e. \textnormal{ implies } \mu \in \mathbf{Reg}.
\end{equation}
This is the so-called Erd\H{o}s-Tur\'an condition. For more details, see \( \cite{STOP} \), Theorem 4.1.1. Finally, we review the concept of Green's functions.
\begin{definition}
	Let \( D \) be a domain in the extended complex plane \( \mathbb{C} \cup \{ \infty \}\), and suppose that \( D \) contains a neighbourhood of \( \infty \). The function \( g_D(\cdotp, \infty): D \to (-\infty, \infty] \) is called the Green's function with pole at infinity for \( D \), if \\
(i) \( g_D(\cdotp, \infty) \) is harmonic on \( D \) and bounded outside neighbourhoods of \( \infty \), \\
(ii)  we have
\[	
	g_D(z,\infty) = \log|z| + O(1), \quad z \to \infty
\]
(iii) \( g_D(z,\infty)  \to 0 \) as \( z \to \zeta \in \partial D \setminus E \), where \( E \) is a set of capacity zero.
\end{definition}
Our goal is to prove the following theorem.
\begin{theorem}\label{main_theorem_light}
	Let \( \Gamma \) be a union of finitely many rectifiable Jordan curves lying exterior to each other and let \( \mu \) be a Borel measure supported on \( \Gamma \) such that the Radon-Nikodym derivative \( d\mu/ds_\Gamma \), where \( s_\Gamma \) denotes the arc length measure, is positive \( s_\Gamma \)-almost everywhere. Let \( z_0 \in \Gamma \) be given, and assume that \( d\mu/ds_\Gamma \) has a jump singularity at \( z_0 \), i.e. (\ref{type_one_singularity}) holds for some \( A, B > 0 \), \( A \neq B \). Then
\[
	\lim_{n \to \infty} n \lambda_{n}(\mu, z_0) = \frac{1}{\frac{d \omega_{\Gamma}(z_0)}{d s_{\Gamma}}} \frac{A-B}{\log A - \log B},
\]
where \( \omega_\Gamma \) denotes the equilibrium measure of \( \Gamma \).
\end{theorem}
Note that if \( A = B \), then changing \( A \) to \( A + \varepsilon \) and letting \( \varepsilon \) to \( 0 \), we obtain 
\[
	\lim_{n \to \infty} n \lambda_n(\mu,z_0) = \frac{d\mu(z_0)}{ds_\Gamma} \frac{ds_\Gamma(z_0)}{d\omega_\Gamma},
\]
which was proved in \( \cite{T1} \). Theorem \( \ref{main_theorem_light} \) is a special case of the following theorem, which we will prove.

\begin{theorem}\label{main_theorem}
	Let \( \mu \in \mathbf{Reg} \) be a Borel measure in  \( \mathbb{C} \) such that \( K = \operatorname{supp}(\mu) \) is compact and if \( \Omega \) denotes the unbounded component of \( \mathbb{C} \setminus K \), then \( \operatorname{cap}(K) = \operatorname{cap}(\operatorname{int}(\mathbb{C} \setminus \Omega)) \). Let \( z_0 \in \partial \Omega \) and suppose that there is an open disk \( D \) with center at \( z_0 \) such that \( D \cap K \) is a \( C^2 \) smooth Jordan arc \( J \) where \( \mu \) is absolutely continuous with respect to the arc length measure \( s_{J} \) and \( d\mu = w_{\mu} ds_{J} \) holds. Also suppose that \( w_{\mu} \) has a jump singularity at \( z_0 \), i.e. \( (\ref{type_one_singularity}) \) holds for some \( A, B > 0 \), \( A \neq B \). Then
\begin{equation}\label{main_result}
	\lim_{n \to \infty} n \lambda_{n}(\mu, z_0) =\frac{1}{\frac{d \omega_{K}(z_0)}{d s_{J}}} \frac{A-B}{\log A - \log B} = 2 \pi \Bigg( \frac{\partial g_{\Omega}(z_0,\infty)}{\partial \mathbf{n}} \Bigg)^{-1} \frac{A-B}{\log A - \log B},
\end{equation}
where \( \omega_K \) is the equilibrium measure of \( K \) and \( \frac{\partial g_{\Omega}(z_0,\infty)}{\partial \mathbf{n}} \) denotes the normal derivative of \( g_{\Omega}(\cdotp,\infty) \) with respect to the inner normal of \( \Omega \).
\end{theorem}

In Section 2, we review the tools we need to prove Theorem \( \ref{main_theorem} \), and we carry out the proof in Section 3 .We mainly follow the method of \( \cite{T1} \), although the jump singularity causes major differences in the proof. The proof consists of three main steps. First, we transform the result from \( [-1,1] \) to \( \mathbb{T} = \{ z \in \mathbb{C}: |z| = 1 \} \), and from this, we obtain asymptotic results for measures supported on lemniscates, i.e. level lines of polynomials. Finally, we prove the general result by approximating Jordan curves with lemniscates.

\section{Tools}

\subsection{Approximating curves with lemniscates}

Let \( T_N \) be a polynomial of degree \( N \). With the notation \( \mathbb{T} = \{ z \in \mathbb{C}: |z|=1 \}, \) we call the level set of the polynomial \( T_N \)
\[
	\sigma = T_{N}^{-1}(\mathbb{T}) = \{ z \in \mathbb{C}: |T_N(z)| = 1 \}
\]
a lemniscate. The set \( L = T_{N}^{-1}(\overline{\mathbb{D}}) \), where \( \mathbb{D} = \{ z \in \mathbb{C}: |z| < 1 \} \) is the open unit disk, is called the closed lemniscate domain. Lemniscates have some useful properties, and we will rely mostly on the following theorem proven in \( \cite{NT} \) by B. Nagy and V. Totik.

\begin{theorem}\label{lemniscate}
	Let \( \Gamma \) consist of finitely many Jordan curves lying exterior of each other, let \( K  \) be the union of the domains enclosed by \( \Gamma \), let \( P \in \Gamma \), and assume that in a neighbourhood of \( P \) the curve is \( C^2 \) smooth. Then for every \( \varepsilon > 0 \) there is a lemniscate \( \sigma_P = \sigma \) consisting of Jordan curves such that \( \sigma \) touches \( \Gamma \) at \( P \), \( \sigma \) contains \( K \) in its interior except for the point \( P \), every component of \( \sigma \) contains precisely one component of \( \Gamma \), and if \( L \) denotes the enclosed lemniscate domain, then
\begin{equation}\label{inner_lemniscate}
	\frac{\partial g_{\mathbb{C}\setminus K}(P, \infty)}{\partial \mathbf{n}} \leq \frac{\partial g_{\mathbb{C}\setminus L}(P, \infty)}{\partial \mathbf{n}} + \varepsilon.
\end{equation}
Also, for every \( \varepsilon > 0 \), there exists a lemniscate \( \sigma \) consisting of Jordan curves such that \( \sigma \) touches \( \Gamma \) at \( P \), \( \sigma \) lies strictly inside of \( K \) except for the point \( P \), \( \sigma \) has exactly one component lying inside every component of \( \Gamma \), and
\begin{equation}\label{outer_lemniscate}
	\frac{\partial g_{\mathbb{C}\setminus L}(P, \infty)}{\partial \mathbf{n}} \leq \frac{\partial g_{\mathbb{C}\setminus K}(P, \infty)}{\partial \mathbf{n}} + \varepsilon.
\end{equation}
\end{theorem}

We will need a certain family of lemniscates to behave in a way that the degree of a polynomial, whose level set yields a lemniscate from this family, is bounded by an absolute constant \( N \). The following theorem grants this in some cases.

\begin{lemma}\label{lemniscate2}
	Let \( \Gamma \) consist of finitely many Jordan curves lying exterior to each other and let \( P \in \Gamma \) and assume that for some disk \( D \) with center at \( P \), \( \Gamma \cap D \) is a \( C^2 \) arc. Then for any \( Q \in \Gamma \) lying sufficiently close to \( P \), the lemniscate \( \sigma_Q \), given by either the first or the second part of Theorem \( \ref{lemniscate} \), can be a translated and rotated copy of \( \sigma_P \).
\end{lemma}
The proof can be found in \cite{T1}, Theorem 2.3. For a more detailed discussion on lemniscates and its properties, see \( \cite{NT} \) and Section 2 in \( \cite{T1} \).

\subsection{Fast decreasing polynomials}

To transform the known results from \( [-1,1] \) to general sets, we need polynomials which imitate the behavior of the Dirac delta, but their degree does not increase too fast. For this purpose, we need the following theorem.
\begin{theorem}\label{fast_decreasing_polynomial_theorem}
	Let \( K \) be a compact subset of \( \mathbb{C} \), let \( \Omega \) be the unbounded component of \( \mathbb{C} \setminus K \) and let \( z_0 \in \partial \Omega \) be an arbitrary point. Suppose that there is a disk \( D \subseteq \Omega \) with \( z_0 \in \partial D \). Then there are constants \( C_0, c_0 > 0 \) such that for every \( n \in \mathbb{N} \), there exists a polynomial \( S_{n,z_0,K} \) of degree at most \( n^{109/110} \) with \( S_{n,z_0,K}(z_0) = 1 \), \( |S_{n,z_0,K}(z)| \leq 1 \) and
\begin{equation}\label{decrease_rate}
	|S_{n,z_0,K}(z)| \leq C_0 e^{-c_0 n^{1/110}}, \quad  |z-z_0| \geq n^{-9/10}, \quad z \notin \Omega.
\end{equation}
Also, if \( z_0 \) lies inside a \( C^2 \) boundary open arc \( J \) of \( \partial \Omega \) such that \( J = \partial \Omega \cap O \) with some open set \( O \), then the constants \( C_0, c_0 \) for which \( (\ref{decrease_rate}) \) holds, can be chosen to be independent of \( z_0 \) lying in any closed subarc of \( J \).
\end{theorem}
The theorem and the proof can be found in \cite{T1}, Corollary 4.2. For a more detailed discussion on fast decreasing polynomial constructions, see \( \cite{IT} \) and Section 4 in \( \cite{T1} \).

\subsection{Properties of Green's functions}

Green's functions are needed because they are, in a sense, strongly connected to equilibrium measures. Since, as Theorem \( \ref{main_theorem} \) states, the asymptotic behavior depends on equilibrium measures, we shall leverage this. The following two theorems describe this connection in detail.
\begin{theorem}\label{pullback_equilibrium_measure}
	Let \( \gamma \) be a Jordan curve, let \( K \) be the closure of the domain enclosed by \( \gamma \) and let \( \Phi \) be the conformal map from \( \mathbb{C}_{\infty} \setminus K \) onto the exterior of the unit disk. Then \( \Phi \) extends to a homeomorphism of \( \gamma \) onto the unit circle \( \mathbb{T} := \{ z \in \mathbb{C}: |z| = 1 \} \) and the equilibrium measure \( \omega_K \) is the pullback of the normalized arc length measure on \( \mathbb{T} \) under \( \Phi \), i.e. for every Borel set \( E \), we have \( \omega_K(E) = |\Phi(E)|/2\pi \), where \( |\cdotp| \) denotes the arc length measure.
\end{theorem}
\begin{theorem}\label{green_equilibrium}
	Let \( K \) be a compact set on the plane, \( \Omega \) the unbounded component of \( \mathbb{C} \setminus K \), and \( J \) a \( C^2 \) arc on the boundary of \( \Omega \) such that \( J = \partial \Omega \cap O \) for some open set \( O \). Assume that \( O \setminus J \) has two components, from which one of them is a subset of \( \Omega \), and the other is a subset of \( K \). Then on \( J \), we have
\[
	d \omega_{K}(z) = \frac{1}{2\pi} \frac{\partial g_{\mathbb{C} \setminus K}(z,\infty)}{\partial \mathbf{n}} d s_J(z).
\]
\end{theorem}
These theorems are proved in Section 3 of \( \cite{T1} \). For a more detailed discussion on equilibrium measures and Green's functions, see Section 3 of \( \cite{T1} \), or \( \cite{TR} \).

\subsection{Polynomial inequalities}

Lastly, we review the tools needed for estimating norms of polynomials.

\begin{theorem}[Nikolskii]\label{nikolskii_inequality}
	Let \( 0 < q \leq p \leq \infty \), let \( \{P_n\} \) be a sequence of polynomials on some interval \( [a,b] \) with \( \deg(P_n) \leq n \), and let \( \{T_n\} \) be a sequence of trigonometric polynomials on \( \mathbb{T} = \{ z \in \mathbb{C}: |z| = 1 \} \) with \( \deg{T_n} \leq n \). Then we have
\begin{equation}
	\| P_n \|_{L^p[a,b]} = O(n^{\frac{2}{q}-\frac{2}{p}}) \| P_n \|_{L^q[a,b]}
\end{equation}
and
\begin{equation}
	\| T_n \|_{L^p(\mathbb{T})} = O(n^{\frac{1}{q}-\frac{1}{p}}) \| T_n \|_{L^q(\mathbb{T})},
\end{equation}
with \( O(n^{\frac{2}{q}-\frac{2}{p}}) \) and \( O(n^{\frac{1}{q}-\frac{1}{p}}) \) independent of \( P_n \) and \( T_n \).
\end{theorem}

Although this result should be stated as an inequality, the preceeding form will be the most useful for us. For the precise statement of the theorem and the proof, see Theorem 4.2.6 in \( \cite{DVLZ} \). For \( C^2 \) Jordan curves and arcs, we have the following theorem.
\begin{theorem}\label{nikolskii_arc}
	Let \( 1 \leq p < \infty \), let \( J \) be a \( C^2 \) Jordan curve or arc, and let \( P_n \) be a polynomial of degree \( n \). Then
	\[
		\sup_{z \in J} |P_n(z)| \leq c n^{2/p} \Bigg( \int |P_n|^p ds_J \Bigg)^{1/p}
	\]
	for some constant \( c > 0 \).
\end{theorem}
For more details and proof, see Theorem 6 in \cite{VA} or Corollary 3.10 in \cite{TV}. 
\begin{remark}
	\textnormal{We often use Theorem \ref{nikolskii_arc} for measures supported on Jordan curves behaving like arc length measures on subarcs. That is, let \( \gamma \) be a Jordan curve, and let \( \gamma_1, \dots, \gamma_k \) be disjoint subarcs of \( \gamma \) such that \( \gamma_1 \cup \dots \cup \gamma_k = \gamma \), and let \( a_1, \dots, a_k \) be positive numbers. Then for the measure \( d\mu(z) = \sum_{l=1}^{k} a_l ds_{\gamma_l}(z) \), we have, for example
\begin{equation}\label{nikolskii_piecewise_arc}
	\| P_n \|_\gamma \leq c n \| P_n \|_{L^2(\mu)}.
\end{equation}
}
\end{remark}

\begin{theorem}[Bernstein]\label{bernstein_inequality}
	Let \( K \subset \mathbb{C} \) be a compact set, let \( \Omega \) be the unbounded component of \( \overline{\mathbb{C}}\setminus K \), and also let \( J \subset \partial \Omega \) be a closed \( C^2 \) arc such that \( J = \partial \Omega \cap O \), where \( O \) is an open set. Let \( J_1 \) be a closed subarc of \( J \) not having common endpoints with \( J \). Then for every \( D > 0 \), there is a constant \( C_D \) such that
\begin{equation}
	|P_{n}^{'}(z)| \leq C_D n \|P_n\|_K, \quad \operatorname{dist}(z,J_1) \leq \frac{D}{n}.
\end{equation}
\end{theorem}

For proof, see Section 7 in \( \cite{T1} \).
\section{The main result}

The proof of Theorem \( \ref{main_theorem} \) is executed in three steps. First, we prove the theorem for a measure supported on the unit circle with a pure jump function as its weight function. For our purposes, it will be enough to know the result on the unit circle for this exact measure, because as we shall see, the asymptotic behavior of the Christoffel function at the boundary will depend mostly on local properties of the measure. Then, by applying the method of polynomial inverse images, we transform the result from the unit circle to lemniscates and then to general curves. For a complete discussion on the polynomial inverse image method, see \( \cite{T2} \). For convenience, we define the pure jump functions
\begin{equation}\label{v_def}
	v_{\mathbb{R}}(t) =
	\begin{cases}
		A & \textnormal{if } t \leq 0, \\
		B & \textnormal{if } t > 0
	\end{cases}
\end{equation}
and
\begin{equation}\label{v_circle_def}
	v_{\mathbb{T}}(e^{it}) =
	\begin{cases}
		A & \textnormal{if } -\pi/2 < t < \pi/2, \\
		B & \textnormal{if } \pi/2 \leq t \leq 3\pi/2,
	\end{cases}
\end{equation}
where \( A,B > 0 \) and \( A \neq B \). Before we turn to the proof of our main result, we need two important lemmas which will be used frequently.
\begin{lemma}\label{majorization}
	With the notations and the assumptions of Theorem \ref{main_result}, we have \( \lambda_n(\mu,z_0) \leq C n^{-1} \) for some constant \( C > 0 \).
\begin{proof}
	Define \( \nu := \left.\mu\right|_{K \setminus J_0} + (\sup_{z \in J_0}w_\mu(z)) s_{J_0} \), where \( J_0 \) is a subarc of \( J \) such that \( z_0 \in J_0 \) and the supremum above is finite. It is clear that \( \lambda_n(\mu,z_0) \leq \lambda_n(\nu,z_0) \), but for \( \nu \), the asymptotic behavior is known. From \cite{T1}, Theorem 1.2, it follows that \( \lambda_n(\nu,z_0) = O(n^{-1}) \)
\end{proof}
\end{lemma}
\begin{lemma}\label{subsequence_estimate}
	Let \( \mu \) be a finite Borel measure on the complex plane and let \( \{ n_k \}_{k=1}^{\infty} \subseteq \mathbb{N} \) be an increasing subsequence of \( \mathbb{N} \) such that we have \( n_{k+1}/n_k \to 1 \). Then
	\[
		c \leq \liminf_{k \to \infty} n_k \lambda_{n_k}(\mu, z) \leq \limsup_{k \to \infty} n_k \lambda_{n_k}(\mu, z) \leq C
	\]
	implies
	\[
		c \leq \liminf_{n \to \infty} n \lambda_{n}(\mu, z) \leq \limsup_{n \to \infty} n \lambda_{n}(\mu, z) \leq C.
	\]
\begin{proof}
	For all \( n \in \mathbb{N} \), define \( n_{k(m)} := \max\{ n_k: n_k \leq m \} \). \( \lambda_{n}(\mu,z) \geq \lambda_{n+1}(\mu,z) \) implies that for all \( n \), we have
\[
	n_{k(m)} \lambda_{n_{k(m)}}(\mu,z) \geq \frac{n_{k(m)}}{m} m \lambda_m(\mu,z) \geq \frac{n_{k(m)}}{n_{k(m)+1}} m \lambda_m(\mu,z)
\]
This implies
\[
	\limsup_{m \to \infty} m \lambda_m(\mu,z) = \limsup_{m \to \infty} \frac{n_{k(m)}}{n_{k(m)+1}} m \lambda_m(\mu,z) \leq \limsup_{m \to \infty} \leq n_{k(m)} \lambda_{n_{k(m)}}(\mu,z) = C.
\]
The proof of
\[
	c \leq \liminf_{m \to \infty} m \lambda_{m}(\mu, z)
\]
is similar.
\end{proof}
\end{lemma}

\subsection{Transferring the result from the interval \( [-1,1] \) to the unit circle}\label{unit_circle}

Let \( \mu_{\mathbb{T}} \) be a Borel-measure on the unit circle such that we have \( d\mu_{\mathbb{T}}(e^{it}) = v_{\mathbb{T}}(e^{it}) ds_{\mathbb{T}}(e^{it}) \), where \( s_{\mathbb{T}} \) denotes the arc length measure on the unit circle and \( v_\mathbb{T} \) was defined in (\ref{v_circle_def}). Note that \( d\mu_\mathbb{T}/ds_{\mathbb{T}} > 0 \) and therefore according to (\ref{Erdos-Turan}), we have \( \mu_\mathbb{T} \in \mathbf{Reg} \). We wish to deduce
\begin{equation}\label{main_result_circle}
	\lim_{n \to \infty} n\lambda_n(\mu_{\mathbb{T}},e^{i\pi/2}) = 2\pi \frac{A-B}{\log A - \log B}
\end{equation}
from results known for the interval \( [-1,1] \). For this purpose, we shall transform the measure \( \mu_{\mathbb{T}} \) into a measure \( \mu_{[-1,1]} \) on the interval \( [-1,1] \) via the mapping \( \varphi: e^{it} \mapsto \cos t \). Note that for us to be able to do this, the symmetry of the measure \( \mu_{\mathbb{T}} \) is needed. With this transformation, we have
 \begin{equation}\label{transform1}
	\int_{-\pi}^{\pi} f(\cos t) d\mu_{\mathbb{T}}(e^{it}) = 2 \int_{-1}^{1} f(x) d\mu_{[-1,1]}(x)
\end{equation}
and if \( d\mu_{[-1,1]}(x) = w_{\mu_{[-1,1]}}(x) dx \), then
\begin{equation}\label{transform2}
	v_{\mathbb{T}}(e^{it}) = w_{\mu_{[-1,1]}}(\cos t)\sin t, \quad t \in [0,\pi).
\end{equation}
According to a result of A. Foulqui\'e Moreno, A. Mart\'inez-Finkelshtein and V. L. Sousa, we have
\begin{equation}\label{result_for_interval}
	\lim_{n \to \infty} n\lambda_n(\mu_{[-1,1]},0) = \pi \frac{A-B}{\log A - \log B}.
\end{equation}
This is a special case of \( \cite{MMS} \), Theorem 11. For more details, see the aforementioned article or \( \cite{NeT} \). To obtain \( (\ref{main_result_circle}) \), we shall establish the lower and upper estimates
\[
	2\pi \frac{A-B}{\log A - \log B} \leq \liminf_{n \to \infty} n\lambda_n(\mu_\mathbb{T},e^{i\pi/2}), \quad  \limsup_{n \to \infty} n\lambda_n(\mu_\mathbb{T},e^{i\pi/2}) \leq 2\pi \frac{A-B}{\log A - \log B}.
\]

\textbf{Lower estimate.} Let  \( S_{2n}(e^{it}) \) be the extremal polynomial for \( \lambda_{2n}(\mu_\mathbb{T},e^{i\pi/2}) \). Define
\[
	P_{n}^{*}(e^{it}) = S_{2n}(e^{it})\Bigg( \frac{1+e^{i(t-\pi/2)}}{2} \Bigg)^{2\lfloor \eta n \rfloor} e^{-i(n+\lfloor \eta n \rfloor)(t-\pi/2)}
\]
and \( P_n(\cos t) = P_{n}^{*}(e^{it}) + P_{n}^{*}(e^{-it}) \). Note that \( P_n(\cos t) \) is a polynomial in \( \cos t \), \( \deg(P_n) \leq n + \lfloor \eta n \rfloor \), and \( P_n(0)  = 1 \).  We wish to estimate the integral \(\int_{-\pi}^{\pi} |P_n(\cos t)|^2 d\mu_\mathbb{T}(e^{it})  \). We do this by splitting the domain of integration into several pieces: we carry out our estimations in the intervals \\ \( [-\pi,-\pi/2-\delta], [-\pi/2-\delta,-\pi/2+\delta], [-\pi/2+\delta,\pi/2-\delta], [\pi/2-\delta,\pi/2+\delta], [\pi/2+\delta,\pi] \), where \( \delta \in (0,\pi/2) \) is arbitrary. First we deal with the integral over \( [\pi/2-\delta,\pi/2+\delta] \). We claim that
\begin{equation}\label{P_n_order}
	|P_n(\cos t)|^2 = |P_{n}^{*}(e^{it})|^2 + O(q^n)
\end{equation} 
holds for some \( |q|<1 \) in \( [\pi/2-\delta,\pi/2+\delta] \). To see this, we have
\[
	|P_n(\cos t)|^2 = |P_{n}^{*}(e^{it}) + P_{n}^{*}(e^{-it})|^2 \leq |P_{n}^{*}(e^{it})|^2 + 2|P_{n}^{*}(e^{it})||P_{n}^{*}(e^{-it})| + |P_{n}^{*}(e^{-it})|^2.
\]
It suffices to show that \( |P_{n}^{*}(e^{it})||P_{n}^{*}(e^{-it})| = O(q^n) \) and \( |P_{n}^{*}(e^{-it})|^2 = O(q^n) \). This follows from a regularity argument and the fact that
\begin{equation}\label{fdp_unit_circle}
	\Bigg( \frac{1+e^{i(t-\pi/2)}}{2} \Bigg)^{2\lfloor \eta n \rfloor} = O(q^n), \quad t \in [-\pi,\pi] \setminus \{ \pi/2 \}.
\end{equation}
Indeed, the regularity of \( \mu_\mathbb{T} \) means that if \( \varepsilon > 0 \) is an arbitrary small number, then there is an \( n_0 \in \mathbb{N} \) such that for all \( n \geq n_0 \), we have
\begin{equation}\label{P_regularity}
	\| P_{n}^{*}\|_{\infty} \leq (1+\varepsilon)^n \|P_{n}^{*}\|_2.
\end{equation}
Since we have \(  \|P_{n}^{*}\|_2 \leq \| S_{2n} \|_2 = O(n^{-1/2}) \), if \( \varepsilon > 0 \) is small enough, it can be seen from (\ref{fdp_unit_circle}) that \( |P_{n}^{*}(e^{it})||P_{n}^{*}(e^{-it})| = O(q^n) \) and  \( |P_{n}^{*}(e^{-it})|^2 = O(q^n) \). Therefore we have
\begin{align*}
	\int_{\pi/2-\delta}^{\pi/2 + \delta} |P_n(\cos t)|^2 d\mu_{\mathbb{T}}(e^{it}) & \leq O(q^n) + \int_{\pi/2-\delta}^{\pi/2 + \delta} |S_{2n}(e^{it})|^2 d\mu_{\mathbb{T}}(e^{it}) \\ & \leq O(q^n) + \lambda_{2n}(\mu_\mathbb{T},e^{i\pi/2}).
 \end{align*}
For the integral \( \int_{-\pi/2-\delta}^{-\pi/2 + \delta} |P_n(\cos t)|^2 d\mu_{\mathbb{T}}(e^{it}) \), using the same argument, we have
\[
	\int_{-\pi/2-\delta}^{-\pi/2 + \delta} |P_n(\cos t)|^2 d\mu_{\mathbb{T}}(e^{it}) \leq O(q^n) + \lambda_{2n}(\mu_\mathbb{T},e^{i\pi/2}).
\]
Now we deal with the integrals over \( [-\pi,-\pi/2-\delta], [-\pi/2+\delta,\pi/2-\delta], [\pi/2+\delta,\pi] \). These three cases are very similar. We claim that \( |P_n(\cos t)| \) is exponentially small in these intervals. Indeed, this can be seen by noticing that \( \big|(1+e^{i(t-\pi/2)})/2\big|^{2\lfloor \eta n \rfloor} \) is also exponentially small and applying the regularity argument \( (\ref{P_regularity}) \) for \( S_{2n} \). Therefore we have
\begin{align*}
	& \int_{-\pi}^{-\pi/2-\delta}|P_n(\cos t)|^2 d\mu_{\mathbb{T}}(e^{it}) = O(q^n), \\
	& \int_{-\pi/2+\delta}^{\pi/2-\delta}|P_n(\cos t)|^2 d\mu_{\mathbb{T}}(e^{it}) = O(q^n), \\
	& \int_{\pi/2+\delta}^{\pi}|P_n(\cos t)|^2 d\mu_{\mathbb{T}}(e^{it}) = O(q^n).
\end{align*}
Combining these estimates and using (\ref{transform1}), we have
\begin{align*}
	\int_{-1}^{1} |P_n(x)|^2 d\mu_{[-1,1]}(x) & = \frac{1}{2} \int_{-\pi}^{\pi} |P_n(\cos t)|^2 d\mu_{\mathbb{T}}(e^{it}) \\
	& \leq O(q^n) + \lambda_{2n}(\mu_\mathbb{T},e^{i\pi/2}).
\end{align*}
Therefore
\begin{align*}
	\liminf_{n \to \infty} 2\Big(1+\frac{\lfloor \eta n \rfloor}{n}\Big)n \lambda_{n+\lfloor \eta n \rfloor}(\mu_{[-1,1]},0) & \leq \liminf_{n \to \infty} \Big(1+\frac{\lfloor \eta n \rfloor}{n}\Big)2n \int_{-1}^{1} |P_n(x)|^2 d\mu_{[-1,1]}(x) \\
	& \leq \liminf_{n \to \infty} \Big(1+\frac{\lfloor \eta n \rfloor}{n}\Big)\Big(O(q^n) + 2n\lambda_{2n}(\mu_\mathbb{T},e^{i\pi/2})\Big)  \\
	& = (1+\eta) \liminf_{n \to \infty} 2n \lambda_{2n}(\mu_\mathbb{T},e^{i\pi/2}).
\end{align*}
This estimate only holds on a certain subsequence of \( \lambda_n(\mu_\mathbb{T}, e^{i\pi/2}) \), but Lemma \( \ref{subsequence_estimate} \) grants us our estimate over the whole sequence. Now applying (\ref{result_for_interval}) yields
\[
	2\pi \frac{A-B}{\log A - \log B} \leq (1+\eta) \liminf_{n \to \infty} n \lambda_{n}(\mu_\mathbb{T},e^{i\pi/2}).
\]
Since \( \eta > 0 \) is arbitrary, we have
\begin{equation}\label{circle_lower_estimate}
	2\pi \frac{A-B}{\log A - \log B} \leq \liminf_{n \to \infty} n \lambda_n(\mu_{\mathbb{T}},e^{i\pi/2}),
\end{equation}
and this is our desired estimate. \\

\textbf{Upper estimate.}
Let \( P_n \) be the extremal polynomial for \( \lambda_n(\mu_{[-1,1]},0) \). Define \( S_n(e^{it}) \) with
\[
	S_n(e^{it}) := P_n(\cos t)\Bigg( \frac{1+e^{i(t-\pi/2)}}{2} \Bigg)^{\lfloor \eta n \rfloor} e^{in(t-\pi/2)},
\]
where \( \eta > 0 \) is an arbitrary small number. Note that \( S_n(e^{it}) \) is a polynomial of degree \( (2+\frac{\lfloor \eta n \rfloor}{n})n \) with \( S_n(e^{i\pi/2}) = 1 \). We wish to give an upper estimate of the integral \( \int_{-\pi}^{\pi} |S_n(e^{it})|^2 d\mu_{\mathbb{T}}(e^{it}) \). Again we do this by splitting domain of integration into several parts, and deal with them separately. First, we estimate the integral over \( [\pi/2-\delta, \pi/2+\delta] \). With \( (\ref{transform1}) \) and \( (\ref{transform2}) \), we have
\begin{align*}
	\int_{\pi/2 - \delta}^{\pi/2 + \delta} |S_n(e^{it})|^2 d\mu_{\mathbb{T}}(e^{it}) & \leq \int_{\pi/2 - \delta}^{\pi/2 + \delta} |P_n(\cos t)|^2 d\mu_{\mathbb{T}}(e^{it}) \\ 
	& = \int_{\cos(\pi/2 + \delta)}^{\cos(\pi/2 - \delta)} |P_n(x)|^2 w_{\mu_{[-1,1]}}(x) \sqrt{1-x^2} dx \\ 
	& \leq \int_{\cos(\pi/2 + \delta)}^{\cos(\pi/2 - \delta)} |P_n(x)|^2 w_{\mu_{[-1,1]}}(x) dx \\ 
	& \leq \lambda_n(\mu_{[-1,1]},0),
\end{align*}
thus we have
\begin{equation}\label{circle_upper_est_1}
	\int_{\pi/2 - \delta}^{\pi/2 + \delta} |S_n(e^{it})|^2 d\mu_{\mathbb{T}}(e^{it}) \leq \lambda_n(\mu_{[-1,1]},0).
\end{equation}
Now we estimate the integral over the intervals \( [-\pi,\pi/2 - \delta] \) and \( [\pi/2+\delta,\pi] \). For this, notice that we have
\begin{equation}
	\Bigg| \frac{1+e^{i(t-\pi/2)}}{2} \Bigg|^{\lfloor \eta n \rfloor} = O(q^n), \quad q < 1, \quad t \in [-\pi,\pi/2-\delta] \cup [\pi/2+\delta,\pi],
\end{equation}
Using (\ref{P_regularity}) similarly as in the lower estimate, we have \( \|P_n\|_{\infty} = O(n^{1/2}) \), therefore
\begin{equation}\label{circle_upper_est_2}
	\int_{-\pi}^{\pi/2-\delta} |S_n(e^{it})|^2 d\mu_{\mathbb{T}}(e^{it}) + \int_{\pi/2+\delta}^{\pi} |S_n(e^{it})|^2 d\mu_{\mathbb{T}}(e^{it}) = O(q^n), \quad |q| < 1. 
\end{equation}
Combining \( (\ref{circle_upper_est_1}) \) and \( (\ref{circle_upper_est_2}) \), we have
\[
	\lambda_{2n + \lfloor \eta n \rfloor}(\mu_{\mathbb{T}}, e^{i\pi/2}) \leq \int_{-\pi}^{\pi} |S_n(e^{it})|^2 d\mu(e^{it}) \leq \lambda_n(\mu_{[-1,1]},0) + O(q^n),
\]
thus by (\ref{result_for_interval}), we obtain
\begin{align*}
	\limsup_{n \to \infty} \Big(2 + \frac{\lfloor \eta n \rfloor}{n}\Big)n \lambda_{2n + \lfloor \eta n \rfloor}(\mu_{\mathbb{T}}, e^{i\pi/2}) & \leq \limsup_{n \to \infty} \Big(2 + \frac{\lfloor \eta n \rfloor}{n}\Big)n\big(O(q^n) + \lambda_n(\mu_{[-1,1]},0)\big) \\ & = (2+\eta) \pi \frac{A-B}{\log A - \log B}.
\end{align*}
Since \( \eta \) is arbitrary, we have
\begin{equation}\label{circle_upper_estimate}
	\limsup_{n \to \infty} n \lambda_n(\mu_{\mathbb{T}}, e^{i\pi/2}) \leq 2\pi \frac{A-B}{\log A - \log B}.
\end{equation}
Now \( (\ref{main_result_circle}) \) follows from \( (\ref{circle_lower_estimate}) \) and \( (\ref{circle_upper_estimate}) \).

\subsection{From the unit circle to lemniscates}\label{lemniscate_case}

Let \( \sigma = \{ z \in \mathbb{C}: |T_N(z)| = 1\} \) be the level line of a polynomial \( T_N \), and assume that \( \sigma \) has no self-intersection. Suppose that \( \deg(T_N) = N \), and let \( z_0 \in \sigma \) be arbitrary. Without the loss of generality, we can assume that \( T_N(z_0) = e^{i\pi/2} \). The points \(T_{N}^{-1}(e^{-i\pi/2}) \) divide \( \sigma \) into \( N \) disjoint arcs \( \sigma_1 \cup \dots \cup \sigma_N \). We can also assume that \( z_0 \in \sigma_1 \). Denote the restriction of \( T_N \) to \( \sigma_i \) with \( T_{N,i} \). Note that except for the possible endpoints, \( T_{N,i}: \sigma_i \to \mathbb{T} \) is one-to-one and therefore invertible. With the definition \( h_{i,j}: \sigma_i \to \sigma_j \), \( h_{i,j}(z) = T_{N,j}^{-1}(T_{N}(z)) \), we introduce the notation \( z_{j} = h_{i,j}(z_{i}) \) for every possible \( i,j \). If \( s_\sigma \) denotes the arc length measure, we have
\begin{equation}
	ds_\sigma(z_j) = ds_\sigma(h_{i,j}(z_i)) = \Bigg| \frac{T_{N}^{'}(z_i)}{T_{N}^{'}(h_{i,j}(z_i))} \Bigg| ds_\sigma(z_i) = \Bigg| \frac{T_{N}^{'}(z_i)}{T_{N}^{'}(z_j)} \Bigg| ds_\sigma(z_i).
\end{equation}
Therefore for all \( f: \sigma \to \mathbb{C} \),
\begin{align*}
	\int_{\sigma_j} f(z_i) |T_{N}^{'}(z_j)| ds_\sigma(z_j) & = \int_{\sigma_i} f(z_i) |T_{N}^{'}(h_{i,j}(z_i))|\Bigg| \frac{T_{N}^{'}(z_i)}{T_{N}^{'}(h_{i,j}(z_i))} \Bigg| ds_\sigma(z_i) \\
	& = \int_{\sigma_i}f(z_i)|T_{N}^{'}(z_i)| ds_\sigma(z_i)
\end{align*}
holds. From this, it follows that
\begin{equation}\label{lemniscate_integral_1}
	\int_{\sigma_j} \Big( \sum_{i=1}^{N} f(z_i) \Big) |T_{N}^{'}(z_j)| ds_\sigma(z_j) = \int_\sigma f(z)|T_{N}^{'}(z)| ds_\sigma(z)
\end{equation}
and
\begin{equation}\label{lemniscate_integral_2}
	\int_{\sigma} \Big( \sum_{i=1}^{N} f(z_i) \Big) |T_{N}^{'}(z_j)| ds_\sigma(z) = N \int_\sigma f(z)|T_{N}^{'}(z)| ds_\sigma(z).
\end{equation}
If \( g:\mathbb{T} \to \mathbb{C} \) is arbitrary, we have
\begin{equation}\label{lemniscate_integral_3}
	\int_\sigma g(T_N(z)) |T_{N}^{'}(z)| ds_\sigma(t) = N \int_{0}^{2\pi} g(e^{it}) dt.
\end{equation}
For more details about lemniscates and its properties, see Section 2 of \( \cite{T1} \).  Another important property of lemniscates is that we can explicitly construct the Green's function with pole at infinity of outer lemniscate domains. To be more precise, we have the following theorem.
\begin{theorem}\label{green_lemniscate}
	Let \( \sigma = \{ z: |T_N(z)| = 1\} \) be the level line of the polynomial \( T_N \) of degree \( N \) and denote the enclosed lemniscate domain with \( L \). Then \\
(i) \( g_{\mathbb{C}\setminus L}(z,\infty) =  \frac{1}{N} \log |T_N(z)| \), \\
(ii) \( \frac{\partial g_{\mathbb{C}\setminus L}(z,\infty)}{\partial \mathbf{n}} = \frac{1}{N} |T_{N}^{'}(z)| \), \\
where \( \frac{\partial g_{\mathbb{C}\setminus L}(z,\infty)}{\partial \mathbf{n}} \) denotes the derivative of \( g_{\mathbb{C}\setminus L}(\cdotp, \infty) \) with respect to the outer normal to \( L \).
\end{theorem}
For more details and proof, see Section 3 of \( \cite{T1} \). We define the measure
\begin{equation}\label{mu_sigma}
d\mu_\sigma(z) = v_{\mathbb{T}}(T_N(z)) ds_\sigma(z),
\end{equation}
where \( s_\sigma \) denotes the arc length measure of \( \sigma \), and \( v_\mathbb{T} \) was defined with (\ref{v_circle_def}). Our goal is to prove
\begin{equation}\label{lemniscate_asymptotic_result}
	\lim_{n \to \infty} n\lambda_n(\mu_\sigma,z_0) = 2\pi\frac{N}{|T_{N}^{'}(z_0)|} \frac{A - B}{\log A - \log B} = 2\pi \Bigg( \frac{\partial g_{\mathbb{C} \setminus L}(z_0,\infty)}{\partial \mathbf{n}} \Bigg)^{-1} \frac{A-B}{\log A - \log B}.
\end{equation}
Note that with the application of Theorem \( \ref{green_equilibrium} \) and Theorem \( \ref{green_lemniscate} \), we have
\[
	\frac{1}{2\pi}\frac{|T_{N}^{'}(z)|}{N} d s_\sigma(z) = d \omega_{L}(z),
\]
therefore \( (\ref{lemniscate_asymptotic_result}) \) can be written as
\[
	\lim_{n \to \infty} n\lambda_n(\mu_\sigma,z_0) = \frac{ds_\sigma(z)}{d\omega_L}  \frac{A - B}{\log A - \log B},
\]
where \( L \) denotes the enclosed lemniscate domain. \\

\textbf{Lower estimate.} Let \( P_n \) be the extremal polynomial for \( \lambda_n(\mu_\sigma,z_0) \), where \( \mu_\sigma \) was defined in (\ref{mu_sigma}), and let \( S_{n,z_0,L}  \) be the fast decreasing polynomial given by Theorem \( \ref{fast_decreasing_polynomial_theorem} \), where \( L \) denotes the lemniscate domain enclosed by \( T_N \). Then \(R_n(z) = P_n(z) S_{n,z_0,L}(z) \) is a polynomial of degree at most \( n + n^{109/110} \), and \( |R_n(z_0)| = 1 \) holds. Let \( \eta > 0 \) be an arbitrary small number, and choose \( \delta > 0 \) such that for every \( z \in \sigma \) with \( |z-z_0| < \delta \), \( z \) lies inside the subarc \( \sigma_1 \), and we have
\begin{equation}\label{T_epsilon_delta_1}
	|T_{N}^{'}(z)| \leq (1+\eta) |T_{N}^{'}(z_0)|.
\end{equation}
Outside of \( |z-z_0| < \delta \), we also have \( |R_n(z)| = o(q^{n^{1/220}}) \) for some \( |q| < 1 \). Indeed, the Nikolskii-type inequality in (\ref{nikolskii_piecewise_arc}) gives us \( \|P_n\|_{\sigma} = O(n^{1/2}) \), and since \( |S_{n,z_0,L}(z)| = O(q^{n^{1/110}}) \) for \( |z-z_0| > \delta \) if \( n \) is large, we obtain that \( |R_n(z)| = o(q^{n^{1/220}}) \) outside of \( |z-z_0| \leq \delta \). Now, the expression \(  \sum_{i=1}^{N} R_n(z_i) \), where  \(\{ z_1, \dots, z_{N}\} = T_{N}^{-1}(T_N(z)) \), is symmetric in the \( z_i \)-s, therefore it is a sum of elementary symmetric polynomials. It follows that there is a polynomial \( Q_n \) of degree at most \( (n+n^{109/110})/N \) such that
\begin{equation}\label{Q_n_def}
	Q_n(T_N(z)) = \sum_{i=1}^{N} R_n(z_i), \quad z \in \sigma.
\end{equation}
For more details on this idea, see Section 5.2 in \( \cite{T1} \) or Section 5 in \( \cite{T2} \).
With the observations made previously, we have
\begin{equation}\label{Q_square}
	|Q_n(T_N(z))|^{2} \leq \sum_{i=1}^{N}|R_n(z_i)|^2 + o(q^{n^{1/220}})
\end{equation}
if the \( \delta > 0 \) choosen previously is small enough, since for any nonequal \( i, j \), the distance \( |z_i - z_j| \) cannot be arbitrarily small.
It is important to note that
\begin{equation}\label{Q_n(T_N(z_0))}
	Q_n(T_N(z_0)) = Q_n(e^{i\pi/2}) = 1 + o(1)
\end{equation}
Since \( |R_n(z)| = o(q^{n^{1/220}}) \) outside of \( |z-z_0| \leq \delta \), applying \( (\ref{decrease_rate}), (\ref{lemniscate_integral_2}) \), \( (\ref{Q_square}) \) yields
\begin{align*}
	\int_{\sigma} |Q_n(T_N(z))|^{2} & |T_{N}^{'}(z)|v_{\mathbb{T}}(T_N(z)) ds_\sigma(z) \\
	& \leq o(q^{n^{1/220}}) + N \int_{|z-z_0| < \delta} |R_n(z)|^2|T_{N}^{'}(z)| v_{\mathbb{T}} (T_N(z)) ds_\sigma(z) \\
	& \leq o(q^{n^{1/220}}) + (1+\eta)N|T_{N}^{'}(z_0)| \lambda_n(\mu_\sigma,z_0).
\end{align*}
On the other hand, with \( (\ref{Q_n(T_N(z_0))}) \), we have
\begin{align*}
	\int_{\sigma} |Q_n(T_N(z))|^{2} & |T_{N}^{'}(z)| v_{\mathbb{T}}(T_N(z)) ds_\sigma(z) \\
	& = N \int_{0}^{2\pi} |Q_n(e^{it})|^2 v_{\mathbb{T}} (e^{it}) dt \\
	& \geq (1+o(1)) N \lambda_{\deg(Q_n)}(\mu_{\mathbb{T}}, e^{i\pi/2}),
\end{align*}
where \( \mu_{\mathbb{T}} \) is defined as in Section \( \ref{unit_circle} \), i.e. \( d\mu_{\mathbb{T}}(e^{it}) = v_{\mathbb{T}}(e^{it}) dt \).
Using \( \deg(Q_n) \leq (n+n^{109/110})/N \), we have
\begin{align*}
	\liminf_{n \to \infty} \deg(Q_n) \lambda_{\deg(Q_n)}(\mu_{\mathbb{T}}, e^{i\pi/2}) & \leq \liminf_{n \to \infty} (1+n^{-1/110}) n /N \lambda_{\deg(Q_n)}(\mu_{\mathbb{T}}, e^{i\pi/2}) \\ & \leq(1+\eta) \liminf_{n \to \infty} \frac{|T_{N}^{'}(z_0)|}{N} n \lambda_n(\mu_\sigma, z_0)
\end{align*}
from which, since \( \eta > 0 \) is arbitrary and \( (\ref{main_result_circle}) \) holds, the application of Lemma \ref{subsequence_estimate} yields the lower estimate
\begin{equation}\label{lemniscate_lower_estimate}
	2\pi \frac{N}{|T_{N}^{'}(z_0)|} \frac{A - B}{\log A - \log B} \leq \liminf_{n \to \infty} n \lambda_n(\mu_\sigma, z_0).
\end{equation}

\textbf{Upper estimate.} Let \( \eta > 0 \) be an arbitrary small number, and similarly to \( (\ref{T_epsilon_delta_1}) \), choose \( \delta > 0 \) such that for every \( z \) with \( |z-z_0| < \delta \),
\begin{equation}\label{T_epsilon_delta_2}
	\frac{|T_{N}^{'}(z_0)|}{1+\eta} \leq |T_{N}^{'}(z)| 
\end{equation}
holds. Let \( Q_n \) be the extremal polynomial for \( \lambda_n(\mu_{\mathbb{T}}, e^{i\pi/2}) \). The Nikolskii-type inequality in (\ref{nikolskii_piecewise_arc}) gives us \( \| Q_n \|_\infty \leq c n \| Q_n \|_{L^2(\mu_\mathbb{T})} \).

Define \( R_n \)  with
\[
	R_n(z) := Q_n(T_N(z)) S_{n,z_0,L}(z),
\]
where \( S_{n,z_0,L}(z) \) is the fast decreasing polynomial given by Theorem \( \ref{fast_decreasing_polynomial_theorem} \) and \( L \) denotes the lemniscate domain enclosed by \( T_N \). Note that \( R_n \) is a polynomial of degree \( \deg(R_n) \leq nN + n^{109/110} \) and we have \( R_n(z_0) = 1 \). Since \( S_{n,z_0,L} \) is fast decreasing and  \( \| Q_n \|_\infty \leq c n \| Q_n \|_{L^2(\mu_\mathbb{T})} = O(n^{1/2}) \), we have
\[
	\sup_{z \in L \setminus \{z:|z-z_0| < \delta \}} |R_n(z)| = o(q^{n^{1/220}}),
\]
for some \( |q| < 1 \). It follows that
\begin{align*}
	\int_{|z_0-z| \geq \delta} |R_n(z)|^2 v_{\mathbb{T}}(T_N(z)) ds_\sigma(z) = o(q^{n^{1/220}}).
\end{align*}
On the other hand, with \( (\ref{T_epsilon_delta_2}) \), we have
\begin{align*}
	\int_{|z_0-z| < \delta} |R_n(z)|^2 v_{\mathbb{T}}(T_N(z))ds_\sigma(z) & \leq \int_{|z_0-z| < \delta} |Q_n(T_N(z))|^2 \frac{|T_{N}^{'}(z)|}{|T_{N}^{'}(z)|} v_{\mathbb{T}}(T_N(z)) ds_\sigma(z) \\
	& \leq \frac{1+\eta}{|T_{N}^{'}(z_0)|} \int_{0}^{2\pi} |Q_n(e^{it})|^2 v_{\mathbb{T}}(e^{it}) dt \\
	& = (1+\eta) \frac{\lambda_n(\mu_{\mathbb{T}}, e^{i\pi/2})}{|T_{N}^{'}(z_0)|}.
\end{align*}
Combining these two observations, we have
\begin{align*}
	\lambda_{\deg(R_n)}(\mu_\sigma, z_0) \leq o(q^n) + (1+\eta) \frac{\lambda_n(\mu_{\mathbb{T}}, e^{i\pi/2})}{|T_{N}^{'}(z_0)|},
\end{align*}
from which
\begin{align*}
	\limsup_{n \to \infty} \deg(R_n) \lambda_{\deg(R_n)}(\mu_\sigma, z_0) & \leq (1+\eta)\limsup_{n \to \infty} \deg(R_n) \frac{\lambda_n(\mu_{\mathbb{T}}, e^{i\pi/2})}{|T_{N}^{'}(z_0)|} \\
	& \leq (1+\eta)\limsup_{n \to \infty} \frac{N}{|T_{N}^{'}(z_0)|} (n + n^{109/110}/N) \lambda_n(\mu_{\mathbb{T}}, e^{i\pi/2})
\end{align*}
follows. Since \( \eta > 0 \) was arbitrary, \( \lambda_n(\mu_{\mathbb{T}},e^{i\pi/2}) = O(n^{-1})\) and we know that
\[
	\lim_{n \to \infty} n \lambda_n(\mu_{\mathbb{T}},e^{i\pi/2}) = 2\pi \frac{A - B}{\log A - \log B}
\]
holds, the application of Lemma \ref{subsequence_estimate} yields the upper estimate
\begin{equation}\label{lemniscate_upper_estimate}
 \limsup_{n \to \infty}  n \lambda_n(\mu_\sigma, z_0) \leq 2\pi \frac{N}{|T_{N}^{'}(z_0)|} \frac{A - B}{\log A - \log B}.
\end{equation}

\subsection{The general case}\label{general}

Now we turn to the proof of the general case. Recall that according to our assumptions in Theorem \( \ref{main_theorem} \), \( \mu \) is a Borel measure in \( \mathbb{C} \) with compact support \(  K = \operatorname{supp}(\mu) \), \( \operatorname{cap}(K) = \operatorname{cap}(\operatorname{int}(\mathbb{C} \setminus \Omega)) \) holds, where \( \Omega \) denotes the unbounded component of \( \mathbb{C} \setminus K \). Also assume that for our \( z_0 \in \partial \Omega \), there is an open disk \( D \) with its center at \( z_0 \) such that \( J := D \cap K \) is a \( C^2 \) smooth Jordan arc. On this \( J \), the measure \( \mu \) is absolutely continuous with respect to the arc length measure \( s_{J} \) and \(  d\mu = w_{\mu} d s_{J} \) holds. Our assumption about the weight function \( w_{\mu} \) is that it has a jump singularity at \( z_0 \) with one-sided limits \( A, B > 0 \), \( A \neq B \), i.e. \( (\ref{type_one_singularity}) \) holds. For convenience, we express \( \mu \) in the form
\begin{equation}
	d\mu(z) = w_\mu(z) ds_J(z) = w_0(z) v_{\mathbb{R}}(f_{J}^{-1}(z)) ds_J(z),
\end{equation}
where \( w_0 \) is continuous, \( w_0(z_0) = 1 \), \( f_J: (-1,1) \to J \) is a one to one parametrization of \( J \) with \( f_J(0) = z_0 \) and \( v_{\mathbb{R}} \) is defined as in \( (\ref{v_def}) \).
Before we turn to our lower and upper estimates, we need a lemma.
\begin{lemma}\label{general_result_lemniscate_lemma}
	With the assumptions and notations made before, for every \( \varepsilon > 0 \), there exists a lemniscate \( \sigma_{z_0} \) such that \( \sigma_{z_0} \) lies strictly inside \( \operatorname{int}(\mathbb{C} \setminus \Omega) \) except for the point \( z_0 \) where \( \sigma_{z_0} \) touches \( J \) and for the enclosed lemniscate domain \( L \) of \( \sigma_{z_0} \) we have
\begin{equation}
	\frac{\partial g_{\mathbb{C} \setminus L}(z_0, \infty)}{\partial \mathbf{n}} \leq \frac{\partial g_{\mathbb{C} \setminus \Omega}(z_0, \infty)}{\partial \mathbf{n}} + \varepsilon.
\end{equation}
In a similar manner, for every \( \varepsilon > 0 \), there exists a lemniscate \( \sigma_{z_0} \) such that \( \operatorname{int}(\mathbb{C} \setminus \Omega) \) is strictly inside the enclosed lemniscate domain \( L \) of \( \sigma_{z_0} \) except for the point \( z_0 \) where \( \sigma_{z_0} \) touches \( J \), and we have
\begin{equation}
	\frac{\partial g_{\mathbb{C} \setminus \Omega}(z_0, \infty)}{\partial \mathbf{n}} \leq \frac{\partial g_{\mathbb{C} \setminus L}(z_0, \infty)}{\partial \mathbf{n}} + \varepsilon.
\end{equation}
\end{lemma}
The proof of this lemma can be found in \( \cite{T1} \), Section 5.3. Note that in the proof, the condition \( \operatorname{cap}(K) = \operatorname{cap(\operatorname{int}(\mathbb{C} \setminus \Omega))} \) was used. \\

\textbf{Lower estimate.} Let \( P_n \) be the extremal polynomial for \( \lambda_n(\mu,z_0) \), and let \( S_{n,z_0,K} \) be the fast decreasing polynomial granted by Theorem \ref{fast_decreasing_polynomial_theorem}. By Lemma \ref{majorization}, we know that \( \| P_n \|_{L^2(\mu)} = O(n^{-1/2}) \). First, we need the following lemma.
\begin{lemma}\label{R_max}
	With the notations and the assumptions of Theorem \ref{main_theorem}, if \( P_n \) is the minimizing polynomial for \( \lambda_n(\mu,z_0) \) and if \( S_{n,z_0,K} \) is the fast decreasing polynomial granted by Theorem \ref{fast_decreasing_polynomial_theorem}, then there exists a polynomial \( H_n \) such that \\
(i) \( H_n(z_0) = 1 \), \\
(ii) \( \deg(H_n) = o(n) \), \\
(iii) \( |H_n(z)| \leq 1 \)  for all \( z \in K \), \\
(iv) \( |P_n H_n S_{n,z_0,K}| \) attains its maximum on \( K \) in the disk \( \Delta_{n^{-9/10}}(z_0) = \{ z \in \mathbb{C}: |z-z_0| \leq n^{-9/10} \} \), and we have
\begin{equation}\label{PH_max}
	\sup_{z \in K \cap \Delta_{n^{-9/10}}(z_0)} |P_n(z)H_n(z)S_{n,z_0,K}(z)| \leq c n^{1/2}
\end{equation}
for some constant \( c \).
\begin{proof}
	The existence of \( H_n \) has been proved in \( \cite{T1} \) Section 6.2, and \( (\ref{PH_max}) \) follows directly from the application of the Nikolskii-type inequality found in Theorem \( \ref{nikolskii_arc} \): since \( d\mu(z) \geq c_0 > 0 \) for every \( z \in J \cap \Delta_{n^{-9/10}} \) for some constant \( c_0 \) if \( n \) is large enough, we have
\begin{equation}\label{nikolskii_for_mu}
	\int_{J \cap \Delta_{n^{-9/10}}(z_0)} |P_n H_n S_{n,z_0,K}|^2 d s_J \leq c_1 \int_{J \cap  \Delta_{n^{-9/10}}(z_0)} |P_n H_n S_{n,z_0,K}|^2 d\mu \leq c_1 \| P_n \|_{L^2(\mu)}
\end{equation}
for some constant \( c_1 \). Therefore, using the maximum modulus principle and Theorem \ref{nikolskii_arc}, we have
\begin{align*}
	\| P_n H_n S_{n,z_0,K} \|_{K} & \leq \| P_n H_n S_{n,z_0,K} \|_{K \cap \Delta_{n^{-9/10}}(z_0)} \leq \| P_n H_n S_{n,z_0,K} \|_{J \cap \Delta_{n^{-9/10}}(z_0)} \\
	& \leq c_2 n \Big( \int_{J \cap \Delta_{n^{-9/10}}(z_0)} |P_n H_n S_{n,z_0,K}|^2 ds_J \Big)^{1/2} \\
	& \leq c_3 n \Big( \int_{J \cap \Delta_{n^{-9/10}}(z_0)} |P_n H_n S_{n,z_0,K}|^2 d\mu \Big)^{1/2} \\
	& \leq c_3 n \| P_n \|_{L^2(\mu)} = O(n^{1/2})
\end{align*}
\end{proof}
\end{lemma}
Note that in Lemma \ref{R_max}, the condition \( \mu \in \mathbf{Reg} \) is used.
Now let \( \sigma_{z_0} \) be a lemniscate inside \( \mathbb{C} \setminus \Omega \) given by the first part of Lemma \( \ref{general_result_lemniscate_lemma} \), and suppose that \( \sigma_{z_0} = \{ z \in \mathbb{C}: |T_N(z)| = 1 \} \) for some polynomial \( T_N \) of degree \( N \) with \( T_N(z_0) = e^{i\pi/2} \). Lemma \( \ref{general_result_lemniscate_lemma} \) grants us that if \( L \) denotes the lemniscate domain enclosed by \( \sigma_{z_0} \), we have
\begin{equation}\label{sigma_green_function}
	\frac{\partial g_{\mathbb{C} \setminus L}(z_0, \infty)}{\partial \mathbf{n}} \leq \frac{\partial g_{\mathbb{C} \setminus \Omega}(z_0, \infty)}{\partial \mathbf{n}} + \varepsilon.
\end{equation}
Let \( H_n \) be the polynomial given by Lemma \ref{R_max}, and define \( R_n \) as
\[
	R_n(z) = P_n(z)S_{n,z_0,K}^{2}(z)H_n(z).
\]
Then \( R_n \) is a polynomial with \( \deg(R_n) = n + o(n) \), \( R_n(z_0) = 1 \), and according to Lemma \ref{R_max}, \( |R_n(z)/S_{n,z_0,K}(z)| \) attains its maximum in \( \Delta_{n^{-9/10}}(z_0) \).
In a small neighbourhood of \( z_0 \), we need the lemniscate \( \sigma_{z_0} \) to be close to the Jordan arc \( J \) in some sense. This is granted by the following lemma.
\begin{lemma}\label{curves_are_quadratically_close}
	Let \( z \in J \) such that \( |z-z_0| \leq 3 n^{-9/10} \), and let \( z^* \in \sigma_{z_0} \) be the closest point to \( z \) such that \( s_J([z_0,z]) = s_{\sigma_{z_0}}([z_0,z^*]) \) holds. Then the mapping \( q(z) = z^* \) is one to one and \( |q(z) - z| \leq C|z_0 - z|^2 \). Also, we have \( ds_J(z) = ds_{\sigma_{z_0}}(z^*) \), thus with the notation \( I_n = \{ z^* \in \sigma_{z_0}: |z^* - z_0| \leq n^{-9/10} \} \), we have
\begin{equation}\label{equal_integrals}
	\int_{z^* \in I_n} |R_n(z^*)|^2 v_{\mathbb{T}}(T_N(z^*))ds_{\sigma_{z_0}}(z^*) = \int_{z^* \in I_n} |R_n(q(z))|^2 v_{\mathbb{R}}(f_{J}^{-1}(z)) ds_{J}(z),
\end{equation}
where \( v_{\mathbb{R}} \) and \( v_{\mathbb{T}} \) are defined as in \( (\ref{v_def}), (\ref{v_circle_def}) \). Therefore
\begin{equation}\label{integrals_are_close}
	\Bigg| \int_{z^* \in I_n} |R_n(z^*)|^2 v_{\mathbb{T}}(T_N(z^*))ds_{\sigma_{z_0}}(z^*) - \int_{z^* \in I_n} |R_n(z)|^2v_{\mathbb{R}}(f_{J}^{-1}(z))ds_{J}(z) \Bigg| \leq C n^{-5/4}
\end{equation}
for some constant \( C \).
\begin{proof}
	We will only prove \( (\ref{integrals_are_close}) \), since the other statements are trivial. First we note that applying the inequality in Theorem \( \ref{bernstein_inequality} \), we obtain that
\[
	\frac{|R_n(q(z)) - R_n(z)|}{|q(z) - z|} \leq C_0 (n+o(n)) \| R_n \|_J,
\]
if \( z \) is close enough to \( z_0 \). The Nikolskii-type inequality in Theorem \( \ref{nikolskii_arc} \) with a reasoning similar to (\ref{nikolskii_for_mu}) (except for \( P_n \) instead of \( P_n H_n S_{n,z_0,K} \)), tells us that
\[
	\| R_n \|_{J} \leq \| P_n \|_J \leq C_1 n \| P_n \|_{L^2(\mu)} \leq C_2 n^{1/2}.
\]
We also have \( |z - z_0| \leq 2n^{-9/10} \), if \( |z^* - z_0| \leq n^{-9/10} \), therefore we obtain
\begin{equation}\label{r(q)-r}
	\begin{aligned}
	|R_n(q(z)) - R_n(z)| & \leq C_3 n^{3/2} |q(z) - z| \\
	& \leq C_4 n^{3/2} |z - z_0|^2 \\
	& \leq C_5 n^{-3/10}
	\end{aligned}
\end{equation}
for all \( z^* \in I_n \). With the application of \( ( \ref{equal_integrals}) \), we have
\begin{align*}
A & = \Bigg| \int_{z^* \in I_n} |R_n(z^*)|^2 v_{\mathbb{T}}(T_N(z^*))ds_{\sigma_{z_0}}(z^*) - \int_{z^* \in I_n} |R_n(z)|^2v_{\mathbb{R}}(f_{J}^{-1}(z))ds_{J}(z) \Bigg| \\
	& = \Bigg| \int_{z^* \in I_n} \Big( |R_n(q(z))|^2 - |R_n(z)|^2 \Big) v_\mathbb{R}(f_{J}^{-1}(z)) ds_J(z) \Bigg|,
\end{align*}
which we need to estimate. The application of the inequalities of H\"older and Minkowski yields
\begin{equation}\label{A_estimate}
\begin{aligned}
	A \leq & \Bigg( \int_{z^* \in I_n} |R_n(q(z)) - R_n(z)|^{2} v_{\mathbb{R}}(f_{J}^{-1}(z))ds_J(z) \Bigg)^{1/2} \\
	& \phantom{gsdgsgsdgsd} \times \Bigg\{ \Bigg( \int_{z^* \in I_n} |R_n(q(z))|^{2} v_{\mathbb{R}}(f_{J}^{-1}(z))ds_J(z) \Bigg)^{1/2} \\
	& \phantom{dgfsdgsgsdgsdgsdgsdgsdgds} + \Bigg( \int_{z^* \in I_n} |R_n(z)|^{2} v_{\mathbb{R}}(f_{J}^{-1}(z))ds_J(z) \Bigg)^{1/2} \Bigg\}. \\
\end{aligned}
\end{equation}
We estimate these integrals one by one. Since \( J \) is a \( C^2 \) Jordan arc, the \( s_J \) measure of \( I_n \) is at most \( C_6 n^{-9/10} \) for some constant \( C_6 \). Therefore for the first of the integrals above, \( (\ref{r(q)-r}) \) yields
\[
	\Bigg( \int_{z^* \in I_n} |R_n(q(z)) - R_n(z)|^{2} v_{\mathbb{R}}(f_{J}^{-1}(z))ds_J(z) \Bigg)^{1/2} \leq C_7 n^{-15/20}.
\]
Since \( P_n \) is the minimizing polynomial for \( \lambda_n(\mu,z_0) \) and \( R_n = P_n S_{n,z_0,K}^{2} H_n \) with \( |S_{n,z_0,K}(z)| \leq 1 \), \( |H_n(z)| \leq 1 \) for all \( z \in K \), we obtain that for the last integral in \( (\ref{A_estimate}) \), we have
\[
	\Bigg( \int_{z^* \in I_n} |R_n(z)|^{2} v_{\mathbb{R}}(f_{J}^{-1}(z))ds_J(z) \Bigg)^{1/2} \leq C_8 n^{-1/2}.
\]
Finally, we have
\begin{align*}
	\Bigg( \int_{z^* \in I_n} & |R_n(q(z))|^{2} v_{\mathbb{R}}(f_{J}^{-1}(z))ds_J(z) \Bigg)^{1/2} \\
	& = \Bigg( \int_{z^* \in I_n} \Big(|R_n(q(z))|^{2} - |R_n(z)|^2 + |R_n(z)|^2 \Big) v_{\mathbb{R}}(f_{J}^{-1}(z))ds_J(z) \Bigg)^{1/2} \\
	& \leq A^{1/2} + C_8 n^{-1/2}.
\end{align*}
Combining these estimates, we obtain that
\begin{align*}
	A & \leq C_7 n^{-15/20} \Big( A^{1/2} + C_8 n^{-1/2} \Big) \\
	& \leq C_9 \max\{n^{-15/20} A^{1/2}, n^{-25/20}\}
\end{align*}
which implies
\[
	A \leq C_{10} n^{-25/20},
\]
and this is what we wanted to prove.
\end{proof}
\end{lemma}
Since \( w_0 \) is continuous and \( w_0(z_0) = 1 \), we have \( w_0(z) = (1 + o(1))w_0(z_0) = 1 + o(1) \) for all \( z \in I_n \) if \( n \) is large enough. Using this and \( (\ref{integrals_are_close}) \), we obtain 
\begin{equation}\label{purple_elephant_1}
\begin{aligned}
	\int_{z^* \in I_n} & |R_n(z^*)|^2 v_{\mathbb{T}}(T_N(z^*))ds_{\sigma_{z_0}}(z^*) \\
	& \leq \frac{1+o(1)}{w_0(z_0)} \int_{z^* \in I_n} |R_n(z)|^2 w_0(z)v_{\mathbb{R}}(f_{J}^{-1}(z))ds_J(z) + Cn^{-5/4} \\
	& \leq (1+o(1)) \int_{z^* \in I_n} |P_n(z)|^2 w_0(z)v_{\mathbb{R}}(f_{J}^{-1}(z))ds_J(z) + Cn^{-5/4}.
\end{aligned}
\end{equation}
By Lemma \( \ref{R_max} \), we have \( \| P_n H_n S_{n,z_0,K} \|_K \leq c n^{1/2} \). Since \( S_{n,z_0,K} \) is fast-decreasing for \( |z-z_0| \geq n^{-9/10} \), we have
\begin{equation}\label{purple_elephant_2}
\begin{aligned}
	\int_{z^* \in \sigma_{z_0} \setminus I_n} |R_n(z^*)|^2 & v_{\mathbb{T}}(T_N(z^*))ds_{\sigma_{z_0}}(z^*) \\
	&  \leq C_0 e^{-c_0 n^{1/110}} \| P_n H_n S_{n,z_0,K} \|_K s_{\sigma_{z_0}}(\sigma_{z_0}) \\
	& \leq D n^{1/2} e^{-c_0 n^{1/110}}
\end{aligned}
\end{equation}
for some constant \( D \). Combining \( (\ref{purple_elephant_1}) \) and \( (\ref{purple_elephant_2}) \) we obtain
\begin{align*}
	\int_{z^* \in \sigma_{z_0}} & |R_n(z^*)|^2 v_{\mathbb{T}}(T_N(z^*))ds_{\sigma_{z_0}}(z^*) \\
	& \leq (1+o(1)) \int_{J} |P_n(z)|^2 w_0(z)v_{\mathbb{R}}(f_{J}^{-1}(z))ds_J(z) + Cn^{-5/4} + D n^{1/2} e^{-c_0 n^{1/110}} \\
	& \leq (1+o(1)) \lambda_n(\mu,z_0) + C n^{-5/4} +  D n^{1/2} e^{-c_0 n^{1/110}}.
\end{align*}

On the other hand, we also have \( \deg(R_n) = n + o(n) \) and
\[
	\lambda_{\deg(R_n)}(\mu_{\sigma_{z_0}},z_0) \leq \int_{z^* \in \sigma_{z_0}} |R_n(z^*)|^2 d\mu_{\sigma_{z_0}}(z^*),
\]
where \( d\mu_{\sigma_{z_0}}(z) = v_{\mathbb{T}}(T_N(z)) ds_{\sigma_{z_0}}(z) \) is a measure on \( \sigma_{z_0} \). To summarize the above, we have
\[
	\lambda_{\deg(R_n)}(\mu_{\sigma_{z_0}},z_0) \leq (1+o(1)) \lambda_n(\mu,z_0) + C n^{-5/4} + D n^{1/2} e^{-c_0 n^{1/110}}.
\]
From this and \( \lambda_n(\mu,z_0) o(n) = o(1) \) it follows that
\[
	\liminf_{n \to \infty} \deg(R_n) \lambda_{\deg(R_n)}(\mu_{\sigma_{z_0}}, z_0) \leq \liminf_{n \to \infty} n \lambda_n(\mu,z_0).
\]
Therefore \( (\ref{sigma_green_function}) \) and \( (\ref{lemniscate_asymptotic_result}) \) yields that
\[
	2\pi \Bigg( \frac{\partial g_{\mathbb{C} \setminus K}(z_0, \infty)}{\partial \mathbf{n}} + \varepsilon \Bigg)^{-1} \frac{A - B}{\log A - \log B} \leq \liminf_{n \to \infty} n\lambda_n(\mu, z_0),
\]
and since \( \varepsilon > 0 \) was arbitrary, the application of Lemma \ref{subsequence_estimate} yields our desired lower estimate
\begin{equation}\label{general_lower_estimate}
	2\pi \Bigg( \frac{\partial g_{\mathbb{C} \setminus K}(z_0, \infty)}{\partial \mathbf{n}} \Bigg)^{-1} \frac{A - B}{\log A - \log B} \leq \liminf_{n \to \infty} n\lambda_n(\mu, z_0)
\end{equation}
follows. \\

\textbf{Upper estimate.} Let \( \sigma_{z_0} = \{ z \in \mathbb{C}: |T_N(z)|=1 \} \) be a lemniscate given by the second part of Lemma \ref{general_result_lemniscate_lemma}, where \( T_N \) is a polynomial of degree \( N \), such that \( \sigma_{z_0} \) touches \( J \) at \( z_0 \), containts \( K \) in its interior, and if \( L \) denotes the enclosed lemniscate domain, we have
\begin{equation}\label{red_giraffe_1}
	\frac{\partial g_{\mathbb{C} \setminus L}(z_0,\infty)}{\partial \mathbf{n}} \geq \frac{\partial g_{\mathbb{C} \setminus K}(z_0,\infty)}{\partial \mathbf{n}} - \varepsilon.
\end{equation}
Without the loss of generality, we can assume that \( T_N(z_0) = e^{i\pi/2} \). As before (see the beginning of Section \( \ref{lemniscate_case} \)), \( T_{N}^{-1}(e^{i\pi/2}) \) divides \( \sigma_{z_0} \) into \( N \) disjoint arcs \( \sigma_1 \cup \dots \cup \sigma_N \), and we can assume that \( z_0 \in \sigma_1 \). Let \( Q_n(z) \) be the extremal polynomial for \( \lambda_n(\mu_{\mathbb{T}}, e^{i\pi/2}) \) where, \( d\mu_{\mathbb{T}}(e^{it}) = v_{\mathbb{T}}(e^{it}) dt \), with the \( v_{\mathbb{T}} \) defined in \( (\ref{v_circle_def}) \). Define
\[
	R_n(z) := Q_n(T_N(z))S_{n,z_0,L}(z),
\]
where \( S_{n,z_0,L} \) is the fast decreasing polynomial given by Theorem \( \ref{fast_decreasing_polynomial_theorem} \). \( R_n \) is a polynomial of degree at most \( nN + n^{109/110} \) and \( |R_n(z_0)| = 1 \). Let \( I_n := \{ z \in J: |z^*-z_0| \leq n^{-9/10} \} \), as defined in Lemma \( \ref{curves_are_quadratically_close} \). With \( (\ref{lemniscate_integral_1}) \), \( (\ref{lemniscate_integral_3}) \) and \( (\ref{integrals_are_close}) \), we obtain
\begin{align*}
	\int_{z^* \in I_n} |R_n(z)|^2 & w_0(z) v_{\mathbb{R}}(f_{J}^{-1}(z)) ds_J(z) \\
	& = (1+o(1))\underbrace{w_0(z_0)}_{= 1}\int_{z^* \in I_n} |R_n(z)|^2 v_{\mathbb{R}}(f_{J}^{-1}(z)) ds_J(z) \\
	 & \leq \frac{1+o(1)}{|T_{N}^{'}(z_0)|}\int_{z^* \in I_n} |Q_n(T_N(z^*))|^2 |T_{N}^{'}(z^*)| v_{\mathbb{T}}(T_N(z^*)) ds_{\sigma_{z_0}}(z^*) + C n^{-5/4} \\
	& \leq \frac{1+o(1)}{|T_{N}^{'}(z_0)|}\int_{\sigma_1} |Q_n(T_N(z^*))|^2 |T_{N}^{'}(z^*)| v_{\mathbb{T}}(T_N(z^*)) ds_{\sigma_{z_0}}(z^*) + C n^{-5/4} \\
	& = \frac{1+o(1)}{|T_{N}^{'}(z_0)|} \int_{0}^{2\pi} |Q_n(e^{it})|^2 v_{\mathbb{T}}(e^{it}) dt + C n^{-5/4} \\
	& = \frac{1+o(1)}{|T_{N}^{'}(z_0)|} \lambda_n(\mu_{\mathbb{T}}, e^{i\pi/2}) + C n^{-5/4} \\
	& = \frac{1+o(1)}{N}\Bigg( \frac{\partial g_{\mathbb{C} \setminus L}(z_0, \infty)}{\partial \mathbf{n}} \Bigg)^{-1} \lambda_n(\mu_{\mathbb{T}}, e^{i\pi/2}) + C n^{-5/4}.
\end{align*}
To estimate \( \int_{J \setminus I_n} |R_n(z)|^2 v_{\mathbb{R}}(f_{J}^{-1}(z)) w_0(z) ds_J(z) \), we shall note that similarly as before, due to the Nikolskii inequality (see Theorem \( \ref{nikolskii_inequality} \) and the remark after that), \( \| Q_n \|_{\infty} = O(n^{1/2}) \). On the other hand, we have \( |S_{n,z_0,L}(z)| \leq C_0 e^{-c_0 n^{1/110}} \) for all \( z \in \sigma_{z_0} \setminus I_n \). Therefore 
\[
	\int_{J \setminus I_n} |R_n(z)|^2 d\mu(z) = o(q^{n^{1/220}})
\]
for some \( |q| < 1 \). Combining these two estimates, we obtain
\begin{equation}\label{red_giraffe_2}
\begin{aligned}
	\lambda_{\deg(R_n)}(\mu,z_0) & \leq \int |R_n(z)|^2 d\mu(z) \\
	& \leq \frac{1+o(1)}{N}\Bigg( \frac{\partial g_{\mathbb{C} \setminus L}(z_0, \infty)}{\partial \mathbf{n}} \Bigg)^{-1} \lambda_n(\mu_{\mathbb{T}}, e^{i\pi/2}) + C n^{-5/4} + o(q^{n^{1/220}}),
\end{aligned}
\end{equation}
and the application \( (\ref{red_giraffe_1}) \), \( (\ref{red_giraffe_2}) \) yields
\begin{align*}
	\limsup_{n \to \infty} \deg(R_n) \lambda_{\deg(R_n)}(\mu, z_0) & \leq \limsup_{n \to \infty} (nN + n^{109/110}) \lambda_{\deg(R_n)}(\mu, z_0) \\
	& \leq \limsup_{n \to \infty} \Bigg( \frac{\partial g_{\mathbb{C} \setminus L}(z_0, \infty)}{\partial \mathbf{n}} \Bigg)^{-1} n \lambda_n(\mu_{\mathbb{T}}, e^{i\pi/2}) \\
	& \leq \limsup_{n \to \infty} \Bigg( \frac{\partial g_{\mathbb{C} \setminus K}(z_0, \infty)}{\partial \mathbf{n}} - \varepsilon \Bigg)^{-1} n \lambda_n(\mu_{\mathbb{T}}, e^{i\pi/2}) \\
	& = 2\pi \Bigg( \frac{\partial g_{\mathbb{C} \setminus K}(z_0, \infty)}{\partial \mathbf{n}} - \varepsilon \Bigg)^{-1} \frac{A - B}{\log A - \log B}
\end{align*}
Since \( \varepsilon > 0 \) was arbitrary, with the application of Lemma \ref{subsequence_estimate} we obtain the desired estimate
\begin{equation}\label{general_upper_estimate}
	\limsup_{n \to \infty} n \lambda_n(\mu, z_0) \leq 2\pi \Bigg( \frac{\partial g_{\mathbb{C} \setminus K}(z_0, \infty)}{\partial \mathbf{n}} \Bigg)^{-1} \frac{A - B}{\log A - \log B}.
\end{equation}
\\[0.5cm]
The combination of \( (\ref{general_lower_estimate}) \) and \( (\ref{general_upper_estimate}) \) yields Theorem \( \ref{main_result} \).

\end{document}